\documentclass{article}

\usepackage{import}
\usepackage{preambleCommands}	
\usepackage{comment}
\usepackage{authblk}

\title{Inexact inner-outer Golub-Kahan bidiagonalization method: A relaxation strategy}
\author[1]{Vincent Darrigrand \thanks {EU Horizon 2020 Project Energy oriented Center of Excellence: toward exascale for energy (EoCoE-II), Project ID: 824158}}
\author[2]{Andrei Dumitrasc \thanks {Bavarian Academic Center for Central, Eastern and Southeastern Europe (BAYHOST)}}
\author[3]{Carola Kruse}
\author[2]{Ulrich Ruede}

\affil[1]{IRIT, CNRS, Toulouse, France}
\affil[2]{Chair for Computer Science 10 - System Simulation, Friedrich-Alexander-Universität Erlangen-Nürnberg, Erlangen, Germany}
\affil[3]{Cerfacs, Toulouse, France}

\begin{document}
\maketitle

\begin{abstract}
  We study an inexact inner-outer generalized Golub-Kahan algorithm for the solution of saddle-point problems with a two-times-two block structure. In each outer iteration, an inner system has to be solved which in theory has to be done exactly. Whenever the system is getting large, an inner exact solver is, however, no longer efficient or even feasible and iterative methods must be used. We focus this article on a numerical study showing the influence of the accuracy of an inner iterative solution on the accuracy of the solution of the block system. Emphasis is further given on reducing the computational cost, which is defined as the total number of inner iterations. We develop relaxation techniques intended to dynamically change the inner tolerance for each outer iteration to further minimize the total number of inner iterations. We illustrate our findings on a Stokes problem and validate them on a mixed formulation of the Poisson problem.
\end{abstract}

\raggedbottom

\begin{comment}
\address[1]{\orgdiv{IRIT}-\orgname{CNRS}, \orgaddress{\state{Toulouse}, \country{France}}}

\address[2]{\orgdiv{Chair for Computer Science 10 - System Simulation, \orgname{©}}, \orgaddress{\state{Erlangen}, \country{Germany}}}

\address[3]{\orgname{Cerfacs}, \orgaddress{\state{Toulouse}, \country{France}}}

\corres{*Andrei Dumitrasc, Friedrich-Alexander-Universität Erlangen-Nürnberg, \orgaddress{Cauerstraße 11,
91058 Erlangen, \country{Germany}}\\ \email{andrei.dumitrasc@fau.de}}

\fundingInfo{EU Horizon 2020 Project Energy oriented Center of Excellence: toward exascale for energy (EoCoE-II), Project ID: 824158 and the Bavarian Academic Center for Central, Eastern and Southeastern Europe (BAYHOST)}
\end{comment}
%\presentaddress{This is sample for present address text this is sample for present address text}

\begin{comment}
\jnlcitation{\cname{%
\author{V. Darrigrand \thanks orgdiv{IRIT}-\orgname{CNRS}, \orgaddress{\state{Toulouse}, \country{France}}}, 
\author{A. Dumitrasc}, 
\author{C. Kruse}, and 
\author{U. R\"ude}} (\cyear{2022}), 
\ctitle{Inexact inner-outer Golub-Kahan bidiagonalization method: a relaxation strategy}, \cjournal{Numerical Linear Algebra with Applications}, \cvol{2017;00:1--6}.}
\end{comment}
\maketitle

\section{Introduction}

Saddle-point systems can be found in a variety of application fields, such as, for example, mixed finite element methods in fluid dynamics or interior point methods in optimization.
An extensive overview about application fields and solution methods for this kind of problems is presented in the well-known article \cite{bgl_2005} by Benzi, Golub and Liesen. In our following study, we want to focus on an iterative solver based on the Golub-Kahan bidiagonalization: the generalized \ac{GKB} algorithm. This solver is designed for saddle-point systems, and was introduced by Arioli\cite{Ar2013}. 
It belongs to the family of Krylov subspace methods and, as such, relies on specific orthogonality conditions, as we will 
review in more detail in \Cref{sec:GKBtheory}.
Enforcing these orthogonality conditions requires solving an \emph{inner problem}, i.e.~formally computing products with matrix inverses (as described in \Cref{alg:GKB}). In practice, this computation is performed with a linear system solver. For this task, we will explore in this article the use of iterative methods to serve as replacement for direct methods that have been used within \ac{GKB} so far. This is essential for very large problems, such as those coming from a discretized \ac{PDE} in 2D or 3D, when direct solvers may reach their limits. 

Using an inner iterative solver might also be advantageous from another point of view as we motivate in the following. 
The solution of large linear systems is often the bottleneck in scientific computing. 
The computational cost and, consequently, the execution time and/or the energy consumption can become prohibitive. 
For the inner-outer iterative \ac{GKB} solver in turn, the principal and costliest part is the solution of the inner system at each outer iteration.  
One approximate metric to measure the cost of the \ac{GKB} solver is the aggregate sum of the number of inner iterations.
For a given setup, the cost of the \ac{GKB} method can hence be optimized by executing only a minimal number of inner iterations necessary for achieving a prescribed accuracy of the solution. 
To reduce this number, there are two possible steps to be taken into account.
In a first step, for a given application it is often unnecessary to solve the linear system with the highest achievable accuracy. This could be the case, for example, in the solution of a discretized \ac{PDE}, when the discretization already introduces an error. A precise solution of the linear system would not improve the numerical solution with respect to the analytic solution of the \ac{PDE} any further than the discretization allows. Next, we come to the second step which will be the main point of the study in this paper. The solution of the inner linear system in the \ac{GKB} method has to be exact, in theory. If we choose a rather low accuracy for the outer iterative solver, an inner exact solution might, however, no longer be necessary, as long as the inner error does not alter the chosen accuracy of the numerical solution. 
This strategy results in a further reduction of the number of inner iterations, since the inner solver will converge in fewer iterations when a less strict stopping tolerance is used. 

In the following study, we address the case where the inner solver has a 
prescribed stopping tolerance and then how this limited accuracy affects the outer process and the quality of its iterates. 
We will show that, with the appropriate choice of parameters, it is possible to make use of inner iterative solvers without compromising the accuracy of the \ac{GKB} result. 
As it can be seen immediately, the lower the accuracy for the inner solver, the less expensive the \ac{GKB} method will be.
Furthermore, we take advantage of the versatility of iterative methods by  adapting the stopping tolerance of the inner solver dynamically. In other words, we prescribe the tolerance of the inner solver according to some criteria determined at each outer iteration. 
This can lead to 
a reduction of the cost, 
since only a minimal number of inner iterations are executed.
Typically, we will reduce the required accuracy for later instances of the inner solver, since later steps of the outer \ac{GKB}-iteration may contribute less to the overall accuracy.

One particular advantage of our proposed method is its generality. The strategy is independent of other choices which are problem-specific, such as the preconditioner for a Krylov method. 
We perform most of our tests on a relatively small Stokes flow problem, to illustrate the salient features.
We confirm our findings by one final test on a larger case of the mixed Poisson problem, including the use of the augmented Lagrangian method, to demonstrate the use in a realistic scenario.

Our study has a similar context as other works on inexact Krylov methods \cite{bouras2000relaxation,bouras2005inexact}, where these algorithms have been investigated from a numerical perspective.
In these articles, the inexactness originates from a limited accuracy of the matrix-vector multiplication or that of the solution of a local sub-problem. 
Similar to what we have described above, it was found that the inner accuracy can be varied from step to step while still achieving convergence of the outer method.
It was shown experimentally that the initial tolerance should be strict, then relaxed gradually, 
with the change being guided by the latest residual norm. 
Other works complemented the findings with theoretical insights, relevant to several algorithms of the Krylov family \cite{simoncini2003theory, simoncini2005relaxed,van2004inexact}. 
It was noted that, in some cases, unless a problem-dependent constant is included, the outer solver may fail to converge if the accuracy of the inner solution is adapted only based on the residual norm. 
This constant can be computed based on extreme singular values, as shown by Simoncini and Szyld \cite{simoncini2003theory}. 
Another source of inexactness can be the application of a preconditioner via an iterative method. 
Van den Eshof, Sleijpen and van Gijzen considered inexactness in Krylov methods originating both from matrix-vector products and variable preconditioning, using iterative methods from the GMRES family \cite{van2005relaxation}. Similarly to earlier work, their analysis relies on the connection between the residual and the accuracy of the solution to the inner problem.
Since applying the preconditioner has the same effect as a matrix-vector product, the same strategies can be applied to more complex, flexible algorithms, such as those involving variable preconditioning: FGMRES \cite{saad1993flexible}, GMRESR \cite{van1994gmresr}, etc. 
A flexible version of the Golub-Kahan bidiagonalization is employed by Chung and Gazzola to find regularized solutions to a problem of image deblurring \cite{chung2019flexible}. In a more recent paper with the same application, Gazzola and Landman develop inexact Krylov methods as a way to deal with approximate knowledge of $\Am$ and $\Am ^T$ \cite{gazzola2021regularization}.
Erlangga and Nabben construct a framework including nested Krylov solvers. They develop a multilevel approach to shift small eigenvalues, leading to a faster convergence of the linear solver \cite{erlangga2008multilevel}.
In subsequent work related to multilevel Krylov methods, Kehl, Nabben and Szyld apply preconditioning in a flexible way, via an adaptive number of inner iterations \cite{kehl2019adaptive}.
Baumann and van Gijzen analyze solving shifted linear systems and, by applying flexible preconditioning, also develop nested Krylov solvers \cite{baumann2015nested}.
McInnes et al. consider hierarchical and nested Krylov methods with a small number of vector inner products, with the goal of reducing the need for global synchronization in a parallel computing setting \cite{mcinnes2014hierarchical}.

Other than solving linear systems, inexact Krylov methods have been studied when tackling eigenvalue problems, as in the paper by Golub, Zhang and Zha \cite{GolZhaZha2000}. 
Although using different arguments, it was shown that the strategy of increasing the inner tolerance is successful for this kind of problem as well.
Xu and Xue make use of an inexact rational Krylov method to solve nonsymmetric eigenvalue problems and observe that the accuracy of the inner solver (GMRES) can be relaxed in later outer steps, depending on the value of the eigenresidual \cite{xu2022inexact}.
Dax computes the smallest eigenvalues of a matrix via a restarted Krylov solver which includes inexact matrix inversion \cite{dax2019restarted}.

Our paper is structured as follows: in \Cref{sec:GKBtheory}, we review the theory and properties of the \ac{GKB} algorithm;  in \Cref{sec:pbDesc}, we describe the specific problem we chose to use as test case for the numerical experiments;  \Cref{sec:constAcc} is meant to illustrate the interactions between the accuracy of the inner solver and that of the outer one in a numerical test setting; \Cref{sec:pertErrAna} describes the link between the error of the outer solver and the perturbation induced by the use of an iterative inner solver. We describe and test our proposed strategy of using a variable tolerance parameter for the inner solver in \Cref{sec:relaxChoices}. We explore the interaction between the method of the \ac{AL} and our strategy in \Cref{sec:AL}. The final section is devoted to concluding remarks.

\section{Generalized Golub-Kahan algorithm}
\label{sec:GKBtheory}

We are interested in saddle-point problems of the form

\begin{align}\label{eqn:spsW}
\left[
\begin{array}{cc}
\Mm & \Am \\ \Am ^T & \mZ     
\end{array}
\right]
\left[
\begin{array}{c}
\wv \\ \pv     
\end{array}
\right]
=
\left[
\begin{array}{c}
\gvv  \\ \rv     
\end{array}
\right]
\end{align}
with $\Mm\in \mathbb{R}^{m\times m}$ being a symmetric positive definite matrix and $\Am\in \mathbb{R}^{m \times n}$ a full rank constraint matrix. 
The generalized \ac{GKB} algorithm for the solution of a class of saddle-point systems was introduced by Arioli \cite{Ar2013}. To apply it to the system (\ref{eqn:spsW}), we first need to have 
the upper block of the right-hand side to be equal to 0. To this end,  we use the transformation
\begin{align}
\label{eq:iniTransf}
\uv &= \wv - \Mm^{-1}\gvv ,\\
\bv &= \rv - \Am^T \uv.
\end{align} 
The resulting system is
\begin{align}\label{eqn:sps}
\left[
\begin{array}{cc}
\Mm & \Am \\ \Am^T & 0     
\end{array}
\right]
\left[
\begin{array}{c}
\uv \\ \pv     
\end{array}
\right]
=
\left[
\begin{array}{c}
0 \\ \bv     
\end{array}
\right],
\end{align}
which is equivalent to that in \Cref{eqn:spsW}. We can recover the $\wv$ variable as $\wv = \uv + \Mm^{-1}\gvv$. 

Let $\Nm\in \mathbb{R}^{n\times n}$ be a symmetric positive definite matrix. To properly describe the \ac{GKB} algorithm, we need to define the following norms
\begin{equation}
    \normM{\vvv}  = \sqrt{ \vvv ^T \Mm \vvv }; \qquad \normN{\qv}  = \sqrt{  \qv ^T \Nm \qv }; \qquad \normNI{\yv}  = \sqrt{ \yv ^T \Nm^{-1} \yv }.
\end{equation}

Given the right-hand side vector $\bv \in \bR  ^n$, the first step of the bidiagonalization is
\begin{equation}
    \label{eq:iniGKBVec}
    \beta _1 = \normNI{\bv}, \quad \qv _1 = \Nm ^{-1} \bv / \beta _1.
\end{equation}
 After $k$ iterations, the partial bidiagonalization is given by
\begin{equation}
\label{eq:oriGKB}
    \begin{cases}
    \Am \Qm _k = \Mm \Vm _k \Bm _k, &\qquad  \Vm _k ^T \Mm \Vm _k = \Id _k  \\
    \Am ^T \Vm _k = \Nm \Qm _k \Bm ^T _k + \beta _ {k+1} \qv _{k+1} \ev _k^T,  &\qquad \Qm _k ^T \Nm \Qm _k = \Id _k
    \end{cases},
\end{equation}
with the bidiagonal matrix
\begin{equation}
\label{eq:BmatGKB}
    \Bm _k=
    \left[ 
	\begin{matrix}
	\alpha_1 & \beta_2 & 0 & \ldots & 0 \\
    0 &	\alpha_2 & \beta_3 &  \ldots & 0 \\
    \vdots &	\vdots & \vdots & \vdots & \vdots  \\
    0 & \ldots & 0 & \alpha_{k-1} & \beta_k \\
    0 & \ldots & 0 & 0 & \alpha_{k}
	\end{matrix}
	\right]
\end{equation}
and the residual term $\beta _ {k+1} \qv _{k+1} \ev _k^T$.
The columns of $\Vm _k$ are orthonormal vectors with respect to the inner product and norm induced by $\Mm$, while the same holds for $\Qm _k$ and $\Nm$ respectively
\begin{equation}
    \begin{split}
       & \vvv _i ^T \Mm \vvv _j = 0, \forall i \neq j; \qquad \normM{\vvv _k} = 1; \\
       & \qv _i ^T \Nm \qv _j = 0, \forall i \neq j; \qquad \normN{\qv _k} = 1.
    \end{split}
\end{equation}
Prior to the normalization leading to $\vvv _k$ and $\qv _k$, the norms are stored as $\alpha _k$ for $\vvv _k$ and $\beta _k $ for $\qv _k$, as detailed in \cref{alg:GKB}.
Using $\Vm _k$, $\Qm _k$ and the relations in \Cref{eq:oriGKB}, we can transform the system from \Cref{eqn:sps} into a simpler form

\begin{equation}\label{eq:transfSys}
	\left[
	\begin{matrix}
		 \Id  _k& \Bm _k   \\
		\Bm _k  ^T & \mZ 
	\end{matrix}
	\right]
	\left[
	\begin{matrix}
		\zv _k   \\
		\yv  _k
	\end{matrix}
	\right]
	=
	\left[
	\begin{matrix}
		\mZ  \\
		\Qm _k ^T \bv 
	\end{matrix}
	\right]. 
\end{equation}
With the choice for $\qv _1$ given in \Cref{eq:iniGKBVec}, we have that $\Qm _k  ^T \bv = \beta _1 \ev _1 $. The solution components to \Cref{eq:transfSys} are then given by
\begin{equation}
\label{eq:zandy}
    \zv _k= \beta _1  \Bm _k ^{-T} \ev _1; \quad \yv _k= - \Bm _k ^{-1} \zv _k,
\end{equation}
where $ \Bm _k ^{-T}$ is the inverse of $ \Bm _k ^{T}$.
We can build the $k$-th approximate solution to \Cref{eqn:sps} as
\begin{equation}
\label{eq:GKBapx}
    \uv _k = \Vm _k \zv _k; \quad \pv _k = \Qm _k \yv _k.
\end{equation}
In particular, after a number of $k=n$ steps and assuming exact arithmetic, we have $\uv _k = \uv$ and $\pv _k = \pv $, meaning we have found the exact solution to \Cref{eqn:sps}. A proof of why $n$ terms are sufficient to find the exact solution is given in the introductory paper by Arioli \cite{Ar2013}. This corresponds to a  scenario where it is necessary to perform the $n$ iterations, although, for specific problems with particular features, the solution may be found after fewer steps. As $ k \rightarrow n$, the quality of the approximation improves ($\uv _k \rightarrow \uv$ and $\pv _k \rightarrow \pv $), with the bidiagonalization residual $\beta _ {k+1} \qv _{k+1} \ev _k^T$ vanishing for $k=n$. 

Given the structure of $\beta _1 \ev _1$ and $ \Bm ^{T}$, we  find 
\begin{equation}
    \label{eq:zetaDef}
    \zeta _1= \frac{\beta _1}{\alpha _1},
	\quad \zeta _k = \zeta _{k-1} \frac{\beta _k}{\alpha _k},
	\quad \zv _k = 	\left[
	\begin{matrix}
		\zv _{k-1} \\
		\zeta _k 
	\end{matrix}
	\right] 
\end{equation}
 in a recursive manner.
 Then, $\uv _k$ is computed as $\uv_k = \uv _{k-1} + \zeta _k \vvv _k $. In order to obtain a recursive formula for $\pv$ as well, we introduce the vector
\begin{equation}
    \dv _k= \frac{\qv _k - \beta _k \dv _{k-1}}{\alpha _k}, \quad \dv _1 = \frac{\qv _1}{\alpha _1}.
\end{equation}
Finally, the update formulas are 
\begin{equation}
\label{eq:GKBupdate}
    \uv _k= \uv _{k-1} + \zeta _k \vvv _k , \quad \pv _k= \pv _{k-1} - \zeta _k \dv _k .
\end{equation}
At step $k$ of \Cref{alg:GKB}, we have the following error in the energy norm.
\begin{equation}
\label{eq:errExact}
    \begin{split}
        \normM{\ev _k} ^2 &= \normM{ \uv _k - \uv } ^2 = \normM{   \Vm _k  \zv _k   - [ \Vm _k \Vm _{n-k}] 
        \left[
	\begin{matrix}
	{\zv _k} \\
	{\zv _{n-k}} 
	\end{matrix}
	\right] 
	} ^2  \\
	&= \normM{ \Vm _{n-k} \zv _{n-k} 	} ^2 
	 =  \normEu{ \zv _{n-k}} ^2  = \sum_{i=k+1}^{n} \zeta _i ^2
    \end{split}
\end{equation}
In the last line, we have made use of the $\Mm$-orthonormality of the $\Vm$ matrices. If we truncate the sum above to only its first $d$ terms, we get a lower bound on the energy norm of the error. The subscript $d$ stands for \textit{delay}, because we can compute this lower bound corresponding to a given step $k$ only after an additional $d$ steps
\begin{equation}
    \label{eq:lowBnd}
    \xi ^2 _{k,d} = \sum_{i=k+1}^{k+d+1} \zeta _i ^2 < \normM{\ev _k} ^2.
\end{equation}
With this bound for the absolute error, we can devise one for the relative error in \Cref{eq:lowBndRel}, which is then used as stopping criterion in \Cref{alg_line:GKBconvCheck} of \Cref{alg:GKB}. 
\begin{equation}
    \label{eq:lowBndRel}
    \bar \xi ^2 _{k,d} = \frac{ \sum_{i=k-d+1}^{k} \zeta _i ^2  }{ \sum_{i=1}^{k} \zeta _i ^2  } .
\end{equation}
The \ac{GKB} algorithm has the following error minimization property. Let $\cV _k = span \{\vvv _1, ..., \vvv _k\}$ and $\cQ _k = span \{\qv _1, ..., \qv _k\}$. Then, for any arbitrary step $k$, we have that
\begin{equation}
\label{eq:errMinProp}
    \min_{\underset{(\Am ^T \uv _k - \bv ) \perp \cQ _k}{
    \uv _k \in  \cV _k,}
    } \normM{ \uv  - \uv _k }
\end{equation}
is met for $\uv _k$ as computed by \Cref{alg:GKB}.

For brevity and because the \ac{GKB} algorithm features this minimization property for the primal variable, our presentation will focus on the velocity for Stokes problems. The stopping criteria for our proposed algorithmic strategies rely on approximations of the velocity error norm. For all the numerical experiments that we have performed, the pressure error norm is close to that of the velocity (less than an order of magnitude apart). In the cases where we operate on a different subspace, as a result of preconditioning, we find that the pressure error norm is actually smaller than that for the velocity. In the case where the dual variable is equally important as the primal, one can use a monolithic approach, such as applying MINRES to the complete saddle-point system.

The \ac{GKB} (as implemented by \Cref{alg:GKB}) is a nested iterative scheme in which each outer loop involves solving an inner linear system. According to the theory given in the paper by Arioli \cite{Ar2013}, the matrices $\Mm$ and $\Nm$ have to be inverted exactly in each iteration. 
We can choose $ \Nm=\frac{1}{\eta} \Id$, whose inversion reduces to a scalar multiplication. In the following sections, unless otherwise specified, we consider $\eta =1 $.
On the other hand, the matrix $\Mm$ depends on the underlying differential equations or the problem setting in general. 
As long as the matrix $\Mm$ is of moderate size, a robust direct solver can be used. 
For large problems, however, a direct solution might no longer be possible and an iterative solver will be required. 
At this point, we face two problems. 
First, depending on the application, inverting $\Mm$ might be more or less costly. 
Second, to achieve a solution quality close to machine precision, an iterative solver might require a considerable number of iteration steps.

\begin{algorithm}
  \caption{Golub-Kahan bidiagonalization algorithm}
  \label{alg:GKB}
  \begin{algorithmic}[1]
  \Require{$\Mm , \Am , \Nm, \bv$, maxit}
  \State{$\beta_1 = \|\bv\|_{\Nm^{-1}}$;  $\qv_1 = \Nm^{-1} \bv / \beta_1$}
  \State{$\wv = \Mm^{-1} \Am \qv_1$; $\alpha_1 = \|\wv\|_{\Mm}$; $\vvv_1 = \wv / \alpha_1$}
  \State{$\zeta_1 = \beta_1 / \alpha_1$; $\dv_1=\qv_1/ \alpha_1$; $\uv^{(1)} =\zeta_{1} \vvv_{1}$; $\pv^{(1)} = - \zeta_1 \dv_1$;    }
  \State{$\bar \xi  _{1,d} = 1;$  $k = 1;$}
  \While{ \textcolor{black}{ $\bar \xi  _{k,d}  >  $ tolerance } and $k < $ maxit }
  \State{$\gvv = \Nm^{-1} \left( \Am^T \vvv_k - \alpha_k \Nm \qv_k  \right) $; $\beta_{k+1} = \|\gvv\|_{\Nm}$}
  \State{$\qv_{k+1} = \gvv / {\beta_{k+1}}$}
  \State{ \textcolor{black}{  $\wv = \Mm^{-1} \left(  \Am \qv_{k+1} - \beta_{k+1} \Mm \vvv_{k} \right)$; } $\alpha_{k+1} = \|\wv\|_{\Mm}$}
  \label{alg_line:innerGKBproblem}
  \State{$\vvv_{k+1} = \wv / {\alpha_{k+1} }$}
  \State{$\zeta_{k+1} = - \dfrac{\beta_{k+1}}{\alpha_{k+1}} \zeta_k$}
  \State{$\dv_{k+1} = \left( \qv_{k+1} - \beta_{k+1} \dv_k \right) / \alpha_{k+1} $}
  \State{$\uv^{(k+1)} = \uv^{(k)} + \zeta_{k+1} \vvv_{k+1}$; $\pv^{(k+1)} = \pv^{(k)} - \zeta_{k+1} \dv_{k+1}$}
 \State{$k = k + 1$}
  \If{$k>d$}
  \State{\textcolor{black}{ $\bar \xi  _{k,d}  = \sqrt{ \sum_{i=k-d+1}^{k} \zeta _i ^2  / \sum_{i=1}^{k} \zeta _i ^2  } $} }
  \label{alg_line:GKBconvCheck}
  \EndIf{}
  \EndWhile{}
  \Return $\uv^{k+1}, \pv^{k+1}$
  \end{algorithmic}
\end{algorithm}

In \Cref{alg_line:innerGKBproblem} of \Cref{alg:GKB}, we have the application of $\Mm^{-1} $ to a vector, which represents what we call the \textit{inner problem}. 
Typically, this is implemented as a call to a direct solver using the matrix $\Mm $  and the vector $\Am \qv_{k+1} - \beta_{k+1} \Mm \vvv_{k}$ as the right hand side. 
The main contribution of this work is a study of the behavior exhibited by \Cref{alg:GKB} when we replace the direct solver employed in \Cref{alg_line:innerGKBproblem} by an iterative one. 
In particular, for a target accuracy of the final \ac{GKB} iterate, we want to minimize the total number of inner iterations. 

Our choice for the inner solver is the unpreconditioned \ac{CG} algorithm, for its simplicity and relative generality. The strategies we propose in the subsequent sections do not rely on any specific feature of this inner solver, and are meant to be applicable regardless of this choice. 
We are interested in reducing the total number of inner iterations in a relative and general manner. 
This is why we do not take preconditioning for \ac{CG} into account, which is usually problem-dependent. 
We measure the effectiveness of our methods based on the percentage of inner iterations saved when compared against a scenario to be described in more detail in the following sections. 

\section{Problem description}
\label{sec:pbDesc}
As test problem, we will use a 2D Stokes flow in a rectangular channel domain $\Omega = \left[-1, L \right] \times \left[-1, 1 \right]$ given by
\begin{equation}
    \label{eq:contStokes}
     \begin{aligned}
    - \Delta \vec{u}  + \nabla p &= 0 \\
 	\nabla \cdot  \vec{u} &=0,
 \end{aligned}
\end{equation}
More specifically, we will consider the Poiseuille flow problem, i.e. a steady Stokes problem with the exact solution 
\begin{equation}
    \label{eq:contStokesSol}
    \begin{cases}
    u_x = 1- y^2, \\
    u_y = 0, \\
    p = -2x + \text{constant}.
    \end{cases}
\end{equation}
The boundary conditions are given as Dirichlet condition on the inflow $\Gamma_{in}= \left\{-1\right\} \times \left[-1, 1 \right]$ (left boundary)
and no-slip conditions on the top and bottom walls
$\Gamma_{c}= \left[-1, L \right] \times \left\{-1\right\} \cup \left[-1, L \right] \times \left\{1\right\} $. The outflow at the right $\Gamma_{out}= \left\{L \right\} \times \left[-1, 1 \right]$ (right) is 
represented as a Neumann condition
\begin{equation*}
	\begin{split}
		\frac{\partial u_x}{\partial x} - p &= 0 \\
		\frac{\partial u_y}{\partial x}  &= 0.
	\end{split}
\end{equation*}
We use Q2-Q1 Finite Elements as discretization method.  
Our sample matrices are generated by the Incompressible Flow \& Iterative Solver Software (IFISS)\footnote{\url{http://www.cs.umd.edu/~elman/ifiss3.6/index.html}} package  \cite{ers07}, see the book by Elman et al. \cite{elman2014finite} for a more detailed description of this reference Stokes problem.

\begin{figure}
    \centering
\begin{tikzpicture}
\begin{axis} [
    width=\FigWid \textwidth,
	height=\FigHei \textwidth,
    minor tick num=3, 
    grid=both,
    xtick = {-1,0,...,5},
    ytick = {-1,0,...,1},
    xlabel = $x$, ylabel = $y$,
    ticklabel style = {font = \scriptsize},
    colormap/jet,
    colorbar,
	enlargelimits=false,
	axis on top,
	axis equal image
]
    \addplot [forget plot] graphics[xmin=-1,xmax=5,ymin=-1,ymax=1] {IFISSCh5SolRaw.png};

\end{axis}
\end{tikzpicture}
\caption{ Exact solution to the Stokes problem in a channel of length 5. Plotted is the $1-y^2$ function, which represents the $x$ direction velocity, overlaid with the mesh resulting from the domain discretization (Q2-Q1 Finite Elements Method). } 
\label{fig:convIFISSCh5SolMesh} 
\end{figure}
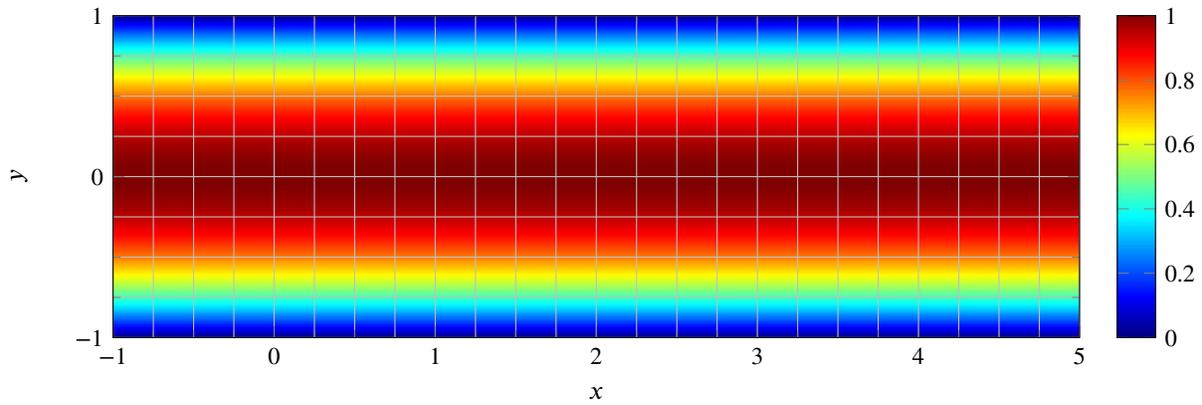 

\begin{figure}
    \centering
    \begin{tikzpicture}
    \begin{axis}[
    ymode=log,
    legend pos= south west,
	table/col sep=tab,
    width=\FigWid \textwidth,
	height=\FigHei \textwidth,
    xlabel={\ac{GKB} iterations},
    ylabel={Velocity error norm}
    ]
    
\addplot+[] table [x=100X, y=100Y]{images/csv/convIFISSChL.csv}; 
\label{fig:item:convIFISSChL_100} 
\addlegendentry{101} 

\addplot+[] table [x=50X, y=50Y]{images/csv/convIFISSChL.csv}; 
\label{fig:item:convIFISSChL_50} 
\addlegendentry{51} 

\addplot+[] table [x=20X, y=20Y]{images/csv/convIFISSChL.csv}; 
\label{fig:item:convIFISSChL_20} 
\addlegendentry{21} 

\addplot+[] table [x=10X, y=10Y]{images/csv/convIFISSChL.csv}; 
\label{fig:item:convIFISSChL_10} 
\addlegendentry{11} 

\end{axis} 
\end{tikzpicture} 
\caption{ \ac{GKB} convergence history for the IFISS channel problem. The length of each channel is given in the legend. Y-axis: Energy norm of the relative error for the velocity.} 
\label{fig:convIFISSChL} 
\end{figure}
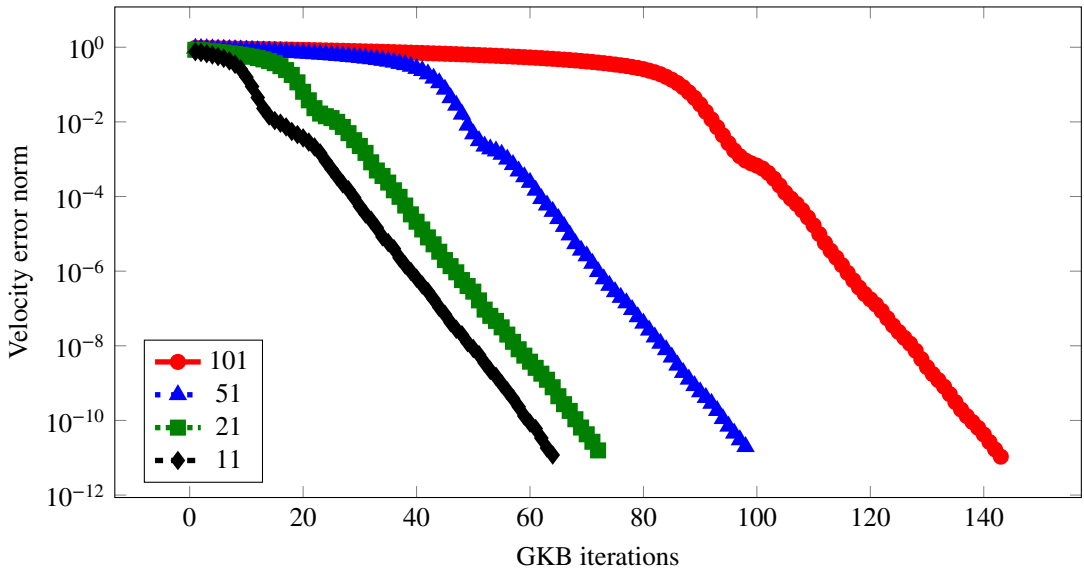

We first illustrate some particular features shown by \ac{GKB} for this problem. We use a direct inner solver here, before discussing the influence of an iterative solver in subsequent sections.
In \Cref{fig:convIFISSChL}, we plot the convergence history for several channels of different lengths, which leads us to noticing the following details. 
The solver starts with a period of slow convergence, visually represented by a plateau, the length of which is proportional to the length of the channel. The rest of the convergence curve corresponds to a period of superlinear convergence, a phenomenon also known for other solvers of the Krylov family, such as \ac{CG}.
The presence of this plateau is especially relevant for our proposed strategies and, since it appears for each channel, we can conclude it is a significant feature of this class of channel problems. In the following numerical examples, we choose as boundary $L=20$ and thus a domain of length 21 units.

\section{Constant accuracy inner solver}
\label{sec:constAcc}

Similar to what has been described by Golub et al. \cite{GolZhaZha2000} for solving eigenvalue problems,
we have observed that when using an iterative method as an inner solver, its accuracy has a clear effect on the overall accuracy of the outer solver (see \Cref{fig:constCGtol}). 

We solve the channel problem described in \Cref{sec:pbDesc} with various configurations for the tolerance of the inner solver, and plot the resulting convergence curves in \Cref{fig:constCGtol}.
The outer solver is always \ac{GKB} with a $10^{-7}$ tolerance. 
The cases we show are: a direct inner solver, three choices of constant inner solver tolerance ($10^{-3}$, $10^{-7}$ and $10^{-8}$), and a final case using a low accuracy solver of ($10^{-3}$) only for the first two iterations, then a high accuracy one ($10^{-14}$).

The stopping criterion for the \ac{GKB} algorithm is a delayed lower bound estimate for the energy norm of the primal variable (see \Cref{eq:lowBnd}). 
As such, \ac{GKB} with a direct inner solver performs a few extra steps, achieving a higher accuracy than the one required, here around $10^{-8}$.

Notice how the outer solver cannot achieve a higher accuracy than that of the inner solver. 
The outer solver stops reducing the error even before reaching the same accuracy as the inner solver.
Replacing the exact inner solver by a \ac{CG} method with a constant tolerance of $10^{-8}$ leads to a convergence process where the error norm eventually reaches a value just below the target accuracy of $10^{-7}$ and does not decrease further. 
This highlights the fact that the inner solver does not need to be exact in order to have \ac{GKB} converge to the required solution. 
For this Poiseuille flow example, however, the inner solver must be least one order of magnitude more precise than the outer one. 

In the last case examined here, we want to see if early imprecise iterations can be compensated later by others having a higher accuracy. 
This strategy of increasing accuracy has been found to work, e.g., in the case of the Newton method for nonlinear problems \cite{dembo1982inexact}.
We tested the case when the first two iterations of \ac{GKB} use an inner solver with tolerance $10^{-3}$, with all the subsequent inner iterations employ a tolerance of $10^{-14}$. 
The resulting curve shows a convergence history rather similar to the case where \ac{CG} has a constant tolerance of $10^{-3}$. 
The outer process cannot reduce the error norm below $10^{-3}$, despite the fact that the bulk of the iterations employ a high-accuracy inner solver. 
This is in correspondence with which was observed by Golub et al. \cite{GolZhaZha2000} for solving eigenvalue problems.

\begin{figure}
    \centering
    \begin{tikzpicture}
    \begin{axis}[
    ymode=log,
    % legend pos= south west,
    legend columns=3,
    legend style={at={(0.5,1.05)}, anchor=south, align=center},
    legend cell align={left},
	table/col sep=tab,
    width=\FigWid \textwidth,
	height=\FigHei \textwidth,
    xlabel={\ac{GKB} iterations},
    ylabel={Velocity error norm},
    extra y ticks={1e-3, 1e-7},
    extra y tick style={draw}
    ]
	\addplot+[] table [x=exactX, y=exactY]{images/csv/constCGtol.csv}; 
\label{fig:item:constCGtol_exact} 
\addlegendentry{direct} 

\addplot+[] table [x=1e-3X, y=1e-3Y]{images/csv/constCGtol.csv}; 
\label{fig:item:constCGtol_1e-3} 
\addlegendentry{\num{1e-3}} 

\addplot+[] table [x=first 2 steps 1e-3 then 1e-14X, y=first 2 steps 1e-3 then 1e-14Y]{images/csv/constCGtol.csv}; 
\label{fig:item:constCGtol_first 2 steps 1e-3 then 1e-14} 
\addlegendentry{mixed\\precision} 

\addplot+[] table [x=1e-7X, y=1e-7Y]{images/csv/constCGtol.csv}; 
\label{fig:item:constCGtol_1e-7} 
\addlegendentry{\num{1e-7}} 

\addplot+[] table [x=1e-8X, y=1e-8Y]{images/csv/constCGtol.csv}; 
\label{fig:item:constCGtol_1e-8} 
\addlegendentry{\num{1e-8}}

\end{axis} 
\end{tikzpicture} 
\caption{\ac{GKB} convergence history for the IFISS Channel test case, depending on the \ac{CG} tolerance (see legend). Y-axis: Energy norm of the relative error for the velocity. Target \ac{GKB} tolerance \num{1e-7}. Mixed precision: first two iterations \num{1e-3}, afterwards \num{1e-14}. The final value for each case with \ac{CG} is: 
\num{3e-3} \ref{fig:item:constCGtol_1e-3}, 
\num{7e-7} \ref{fig:item:constCGtol_1e-7}, 
\num{8e-8} \ref{fig:item:constCGtol_1e-8} , 
\num{2e-3} \ref{fig:item:constCGtol_first 2 steps 1e-3 then 1e-14}.
Only the cases \ref{fig:item:constCGtol_1e-8}  and \ref{fig:item:constCGtol_exact} converge successfully, reducing the error norm below \num{1e-7}.} 
\label{fig:constCGtol} 
\end{figure} 

An interesting observation is that all the curves in \Cref{fig:constCGtol} overlap in their initial iterations, until they start straying from the apparent profile, eventually leveling off. In \Cref{sec:pertErrAna}, we analyze the causes leading to these particular behaviors and link them to the accuracy of the inner solver.

\section{Perturbation and error study}
\label{sec:pertErrAna}
In this section we describe how the error associated with the iterates of \Cref{alg:GKB} behaves if we use an iterative solver for the systems involving $\Mm ^{-1}$. 
We can think of the approximate solutions of these inner systems as perturbed versions of those we would get when using a direct solver. 
The error is then characterized in terms of this perturbation and the implications motivate our algorithmic strategies given in the subsequent sections. With this characterization, we can also explain the results in \Cref{sec:constAcc}.

The use of an iterative inner solver directly affects the columns of the $\Vm$ matrix. In the following, $\Vm$ denotes the unperturbed matrix, with $\Em _{\Vm}$ being the associated perturbation matrix. In particular, we are interested in the $\Mm$ norm of the individual columns of $\Em _{\Vm}$, which gives us an idea of how far we are from the \enquote{ideal} columns of $\Vm$. 

Changes in the  $\vvv$ and $\qv$ vectors also have an impact on their respective norms $\alpha$ and $\beta$, which shift away from the values they would normally have with a direct inner solver. In turn, these changes propagate to the coefficients $\zeta$ used to update the iterates $\uv$ and $\pv$. 
Our observations concern the $\zv$ vector, its perturbation $\ev _{\zv}$ and their effect on the error of the primal variable $\uv$ measured in the $\Mm$ norm. 
The entries of $\zv$ change sign every iteration, but we will only consider them in absolute value, as it is their magnitude which is important. 
In the following, we will denote perturbed quantities with a hat.

\subsection{High initial accuracy followed by relaxation}
\label{subsec:theoPrecRelax}

In this subsection, we take a closer look at the interactions between the perturbation and the error. For us, perturbation is the result of using an inexact inner solver and represents a quantity which can prevent the outer solver from reducing the error below a certain value. The error itself needs to be precisely defined, as it may contain several components, each minimized by a different process. Because we focus on the difference between the perturbed and the unperturbed \ac{GKB}, sources of error that affect both versions, such as the round-off error, are not included in the following discussion. According to the observations by Jir{\'a}nek and Rozlo{\v{z}}n{\'\i}k, the accuracy of the outer solver depends primarily on that of the inner solver, since the perturbations introduced by an iterative solver dominate those related to finite-precision arithmetic \cite{jiranek2008maximum}. We take the exact solution $\uv$ to be equal to $\uv _n$, the $n$-th iterate of the unperturbed \ac{GKB} with exact arithmetic.

At step $k$ of the \ac{GKB}, we have the error, 
\begin{align}
    \normM{\ev _k} &= \normM{ \hat \uv _k - \uv }, 
\end{align}
where $\hat  \uv _k$ is the current approximate solution and $\uv $ is the exact one. Both can be written as linear combinations of columns from $\Vm$ with coefficients from $\zv$. Let $\hat \uv _k$ come from an inexact version of \Cref{alg:GKB}, where the solution of the inner problem (a matrix-vector product with $\Mm ^{-1}$) includes perturbations. The term $ \uv = \Vm _n \zv _n$ is available after $n$ steps of \Cref{alg:GKB} in exact arithmetic, without perturbations. We separate the first $k$ terms, which have been computed, from the remaining ($n-k$).   

\begin{equation}
\label{eq:errPert}
    \begin{split}
        \normM{\ev _k} ^2 &= \normM{ \hat \uv _k - \uv } ^2= \normM{  ( \Vm _k + \Em _{\Vm} ) ( \zv _k + \ev _{\zv} ) - [ \Vm _k \Vm _{n-k}] 
        \left[
	\begin{matrix}
	{\zv _k} \\
	{\zv _{n-k}} 
	\end{matrix}
	\right] 
	} ^2 \\
	&= \normM{ \Em _{\Vm} \zv _k + \Em _{\Vm} \ev _{\zv} + \Vm _k \ev _{\zv} - \Vm _{n-k} \zv _{n-k} 	} ^2 \\
	& \leq \normM{ \Em _{\Vm} \zv _k } ^2 + \normM{ \Em _{\Vm} \ev _{\zv} } ^2 +   \normEu{ \ev _{\zv}  } ^2 + \normEu{ \zv _{n-k}} ^2
    \end{split}
\end{equation}
In the last line, we have made use of the $\Mm$-orthonormality of the $\Vm$ matrices. 

In the case of a direct inner solver, we can leave out the perturbation terms, recovering the result $\normM{\ev _k} ^2 = \normEu{ \zv _{n-k}} ^2 = \sum_{i=k+1}^{n} \zeta _i^2 $ given by Arioli \cite{Ar2013}. This is simply the  error coming from approximating $\uv $ (a linear combination of $n$ $\Mm$-orthogonal vectors) by $ \uv _k $ (a linear combination of only $k$ $\Mm$-orthogonal vectors). This term decreases as we perform more steps of \Cref{alg:GKB} ($ k \rightarrow n $). By truncating the sum $\sum_{i=k+1}^{n} \zeta _i ^2$, we obtain a lower bound for the squared error.

The remaining three terms in \Cref{eq:errPert} include the perturbation coming from the inexact inner solution. Our goal is to minimize the total number of iterations of the inner solver, so we are interested in knowing how large can these terms be allowed to be, such that we still recover a final solution of the required accuracy. The answer is to keep them just below the final value of the fourth one, $\normEu{ \zv _{n-k}}$, below the acceptable algebraic error. If they are larger, the final accuracy will suffer. If they are significantly smaller, then our inner solver is unnecessarily precise and expensive. 

The following observations rely on the behavior of the $\zv$ vector. At each iteration, this vector gains an additional entry, while leaving the previous ones unchanged. These entries form a (mostly) decreasing sequence and have a magnitude below 1 when reaching the superlinear convergence phase. Unfortunately, we cannot yet provide a formal proof of these properties, but having seen them consistently reappear in our numerical experiments encourages us to consider them for motivating our approach. These properties appear in both cases, with and without perturbation. 

The decrease in the entries of the coefficient vector used to build the approximation has also been observed and described for other Krylov methods (see references \cite{van2004inexact,simoncini2003theory,simoncini2005relaxed}). Their context is that of inexact matrix-vector products, which is another way of viewing our case. The fact that new entries of $\zv$ are simply appended to the old ones and that they are smaller than one is linked to the particular construction specific to \ac{GKB}.

Back to \Cref{eq:errPert}, let us assume the perturbation at each iteration is constant, i.e. the $\Mm$ norm of each column of $\Em _{\Vm}$ is equal to the same constant. 
Then, the vector $\Em _{\Vm} \zv _k$ will be a linear combination of perturbation vectors with coefficients from $\zv _k$. Following our observations concerning the entries of $\zv _k$, the first terms of the linear combination will be the dominant ones, with later terms contributing less and less to the sum. If the perturbation of the first $\vvv$ has an $\Mm$ norm below our target accuracy, the term $ \normM{ \Em _{\Vm} \zv _k } $ will never contribute to the error.
We can allow the  $\Mm$ norm of the columns of $\Em _{\Vm}$ to increase, knowing the effect of the perturbation will be reduced by the entries of $\zv$, which are decreasing and less than one. The \ac{GKB} solution can be computed in a less expensive way, as long as the term $\normM{ \Em _{\Vm} \zv _k }$ is kept below our target accuracy. The perturbation should initially be small, then allowed to increase proportionally to the decrease of the entries in $\zv$.

Next, we describe the terms including $\ev _{\zv}$. Let the following define the perturbed entries of $\hat \zv$
\begin{equation*}
        \hat \zeta _k = - \hat \zeta _{k-1} \frac{\hat \beta _k}{\hat \alpha _k} = - \hat \zeta _{k-1} ( \frac{ \beta _k}{ \alpha _k} + \epsilon _k).
\end{equation*}

The term $\epsilon _k$ is the perturbation introduced at iteration $k$, coming from the shifted norms associated with $\qv _k$ and $\vvv _k$. This term is then multiplied by $\hat \zeta _{k-1}$ which, according to our empirical observations, decreases at (almost) every step. If we assume $\epsilon _k$ is constant, the entries of $\ev _{\zv}$ decrease in magnitude and the norm $ \normEu{ \ev _{\zv} }$ is mostly dominated by the first vector entry. The strategy described for the term $ \normM{ \Em _{\Vm} \zv _k } $ also keeps $ \normEu{ \ev _{\zv} }$ small. We start with a perturbation norm below the target accuracy, to ensure the quality of the final iterate. Gradually, we allow an increase in the perturbation norm proportional to the decrease of $\hat \zeta _k$ to reduce the costs of the inner solver. Finally, since the vector $\ev _{\zv} $ decreases similarly to $ \zv $, the term $ \normM{ \Em _{\Vm} \ev _{\zv} }$ can be described in the same way as $\normM{ \Em _{\Vm} \zv _k }$.

We close this section by emphasizing the important role played by the first iterations and how the initial perturbations can affect the accuracy of the solution. Notice that the perturbation terms included refer to all the $k$ steps, not just the latest one. Relaxation strategies that start with a low accuracy and gradually increase it are unlikely to work for \ac{GKB} and other algorithms with similar error minimization properties. Since the first vectors computed are the ones that contribute the most to reducing the error, they should be determined as precisely as possible. Even if we follow a perturbed iteration exclusively by very accurate ones, this will not prevent the perturbation from being transmitted to all the subsequent vectors, and potentially be amplified by multiplication with matrices and floating-point error. With these observations in mind, we can understand the results in \Cref{sec:constAcc}. 

These findings are in line with those concerning other Kylov methods in the presence of inexactness (see Section 11 of the survey by Simoncini and Szyld \cite{simoncini2007recent} and the references therein). \ac{GKB} is not the only method which benefits from lowering the accuracy of the inner process, and the reason why this is possible is linked to the decreasing entries of the coefficient vector.
 
 \section{Relaxation strategy choices}
\label{sec:relaxChoices}

We have seen in \Cref{subsec:theoPrecRelax} that we can allow the perturbation norm to increase in a safe way, as long as the process is guided by the decrease of $ \abs{ \hat \zeta } $. This means that we can adapt the tolerance of the inner solver, such that each call is increasingly cheaper, without compromising the accuracy of the final \ac{GKB} iterate. Then, at step $k$ we can call the inner solver with a tolerance equal to $\tau / f(\zeta)$. The scalar $\tau$ represents a constant chosen as either the target accuracy for the final \ac{GKB} solution, or something stricter, to counteract possible losses coming from floating-point arithmetic. The function $f$ is chosen based on the considerations described below, with the goal of minimizing the number of inner iterations. 

A similar relaxation strategy was used in a numerical study by Bouras and Frayss\'e \cite{bouras2005inexact} to control the % amount
magnitude of the perturbation introduced by performing inexact matrix-vector products. They employ Krylov methods with a residual norm minimization property, so the proposed criterion divides the target accuracy by the latest residual norm. In our case, because of the  minimization property in \Cref{eq:errMinProp}, we need to use the error norm instead of the residual, since it is the only quantity which is strictly decreasing. Due to the actual error norm being unknown, we rely on approximations found via $\zeta$.

Considering the error characterization of the unperturbed process $\normM{\ev _k} ^2 = \sum_{i=k+1}^{n} \zeta _i ^2$, we can approximate the error by the first term of the sum, which is the dominant one. However, when starting iteration $k$ we do not know $\zeta _{k+1}$, not even $\zeta _{k}$, so we cannot choose a tolerance for the inner solver required to compute $\uv _k$ based on these. What we can do is predict these values via extrapolation, using information from the known values $\zeta _{k-1}$ and $\zeta _{k-2}$. We know that in general $ \frac{\beta _k}{\alpha _k} = \frac{\zeta _k}{\zeta _{k-1}} $ acts as a local convergence factor for the $\abs{\zeta }$ sequence. We approximate the one for step $k$ by using the previous one $\frac{\zeta _{k-1}}{\zeta _{k-2}}$. Then, we can compute the prediction $\tilde \zeta _k := \zeta _{k-1} \frac{\zeta _{k-1}}{\zeta _{k-2}}$. By squaring the local convergence factor, we get an approximation for $ \zeta _{k+1} $ as  $\tilde \zeta _{k+1} := \zeta _{k-1} \left(  \frac{\zeta _{k-1}}{\zeta _{k-2}} \right) ^2$, which we can use to approximate $ \normM{\ev _k} $ and adapt the tolerance of the inner solver.

In practice, we only consider processes which include perturbation, and assume we have no knowledge of the unperturbed values $\abs{\zeta }$. As such, for better readability, we drop the hat notation with the implicit convention that we are referring to values which do include perturbation and use them in the extrapolation rule above.

For some isolated iterations, it is possible that $\abs{\zeta _k }  \geq \abs{\zeta _{k-1} }  $. This behavior is then amplified through extrapolation, potentially leading to even larger values. In turn, this can cause an increase in the accuracy of the inner solver, following a stricter value for the tolerance parameter $\tau / f(\zeta)$. In \Cref{subsec:theoPrecRelax}, we have shown that there is no benefit in increasing this accuracy. The new perturbation would be smaller in norm, but the error $\normM{\ev _k}$ would be dominated by the previous, larger perturbation. As such, we propose computing several candidate values for the stopping tolerance of the inner solver, and choose the one with maximum value. Since these are only scalar quantities, the associated computational effort is negligible, but the impact of a well-chosen tolerance sequence can lead to significant savings in the total number of inner iterations. The candidate values are:
\begin{equation}
\label{eq:relaxChoices}
\begin{cases}
	\text{the value at the previous step}, \\
	 \tau / \abs{ \zeta _{k-1} } , \\
	 \tau / \abs{ \tilde \zeta _{k} } , \\
	 \tau / \abs{ \tilde \zeta _{k+1} } .
	\end{cases}
\end{equation}

To prevent a limitless growth of the tolerance parameter, we impose a maximum value of $ 0.1 $. All these choices are safe in the sense that they do not lead to introduction of perturbations which prevent the outer solver from reaching the target accuracy. 

We proceed by testing these relaxations strategies on the problem described in \Cref{sec:pbDesc}.  The initial tolerance for \ac{CG} is set to  $\tau = 10^{-8}$, one order of magnitude more precise than the one set of \ac{GKB}. As a baseline for comparison, we first keep the tolerance constant, equal to $\tau$. Then, we introduce adaptivity using $\tau / \abs{ \zeta_{k-1}  }$. The third case changes the tolerance according to $\tau / \abs{ \tilde \zeta _{k+1} }$, the latter term being a predicted approximation of the current error. Finally, we employ a hybrid approach, where all candidate values in \Cref{eq:relaxChoices} are computed, but only the largest one is used. In the legends of the following plots, these four cases are labeled \ConstantCase, \AdaptiveCase, \PredictedCase, and \HybridCase, respectively. To monitor \ac{GKB} convergence, we track the lower bound for the energy norm of the error corresponding to the primal variable given in \Cref{eq:lowBnd}. For easy reference, all the choices used and their respective labels are given below. We define $\tau = 10^{-8}$.
\begin{align}
    \mathtt{(\ConstantCase)}	&: \tau, \label{eq:cst}\\
	(\mathtt{\AdaptiveCase})  &: \nicefrac{\tau}{\abs{ \zeta _{k-1} }} , \label{eq:z}\\	
	(\mathtt{\PredictedCase}) &: \nicefrac{\tau}{\abs{ \tilde \zeta _{k+1} }} ,\label{eq:adasquare}\\
	(\mathtt{\HybridCase}) &: \max\left\{\nicefrac{\tau}{\abs{ \zeta _{k-1} }}, \nicefrac{\tau}{\abs{ \tilde \zeta _{k} }}, \nicefrac{\tau}{\abs{ \tilde \zeta _{k+1} }},  \text{previous value} \right\} .\label{eq:hybrid}\\
		(\mathtt{\OptimalCase})  &: \nicefrac{\tau}{ ( parameter \cdot \abs{ \zeta _{k-1} } )} , \label{eq:zOptim}
\end{align}

Only the last scenario above, \OptimalCase{}, is left to explain. To see if the parameter-free choices can be improved, we run one more case which includes adaptivity by using $\abs{ \zeta _{k-1}  }$, but also one constant parameter tuned experimentally. This is motivated by the fact that the considerations leading to \Cref{eq:relaxChoices} rely mostly on approximations and inequalities, which means we have an over-estimate of the error. It may be possible to reduce the total number of iterations further, by including an (almost) optimal, problem-dependent constant. The goal is to find a sequence of tolerance parameters with terms that are as large as possible, while guaranteeing the accuracy of the final \ac{GKB} iterate.

\begin{figure}
	\centering
	\begin{tikzpicture}
		\begin{axis}[
			ymode=log,
% 			legend pos= north east,
			legend style={at={(0.5,1.05)}, anchor=south},
			legend columns=3,
			table/col sep=tab,
			width=\FigWid \textwidth,
			height=\FigHei \textwidth,
			xlabel={Inner \ac{CG} iterations},
			ylabel={Lower bound}
			]
			
			\addplot+[name path=A] table [x=cumul iter, y=lower bound]{images/csv/IFISS20/original/const.csv}; 
			\addlegendentry{ \ConstantCase} 
			
			\addplot+[] table [x expr=\thisrow{cumul iter}, y=lower bound]{images/csv/IFISS20/original/hybrid.csv}; 
			\addlegendentry{\HybridCase }

			\addplot+[] table [x expr=\thisrow{cumul iter}, y=lower bound]{images/csv/IFISS20/original/zOptimal.csv}; 
			\addlegendentry{ \OptimalCase} 
			
			\addplot+[] table [x expr=\thisrow{cumul iter}, y=lower bound]{images/csv/IFISS20/original/zAdaSq.csv}; 
			\addlegendentry{ \PredictedCase} 
			
			\addplot+[] table [x expr=\thisrow{cumul iter}, y=lower bound]{images/csv/IFISS20/original/z.csv}; 
			\addlegendentry{ \AdaptiveCase} 
			
			\addplot[name path=B, draw=none] table [x expr=0.7*\thisrow{cumul iter}, y=lower bound]{images/csv/IFISS20/original/const.csv} node[pos=1]{30\%};
			
			\addplot[black!50, opacity=0.5] fill between[of=A and B];
			
			\addplot[name path=C, draw=none] table [x expr=0.6*\thisrow{cumul iter}, y=lower bound]{images/csv/IFISS20/original/const.csv} node[pos=1]{40\%};
			\addplot[black!30, opacity=0.5] fill between[of=A and C];

% 			\addplot[name path=C, draw=none] table [x expr=0.35*\thisrow{cumul iter}, y=lower bound]{images/csv/IFISS20/original/const.csv} node[pos=1]{65\%};
% 			\addplot[black!20, opacity=0.5] fill between[of=A and C];

		\end{axis} 
	\end{tikzpicture} 
	\caption{Lower bound (\Cref{eq:lowBnd}) for the error norm associated with the \ac{GKB} iterates versus the cumulative number of inner \ac{CG} iterations when solving the original problem from \Cref{sec:pbDesc}. The parameter used in \OptimalCase{} is $0.05$. See \Cref{eq:cst,eq:z,eq:zOptim,eq:adasquare,eq:hybrid} for the strategies denoted by the labels.}
	\label{fig:originallower bound}
\end{figure}
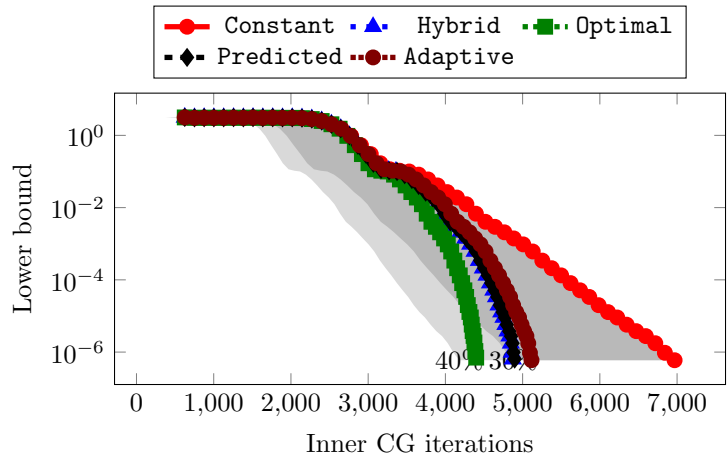

All the results are given in \Cref{tab:varTolZ} and \Cref{fig:originallower bound}. \HybridCase{} offers the highest savings among the parameter-free choices (30\%), but \OptimalCase{}, the test with the problem-dependent constant, reveals that we can still improve this performance by about 6\%.

\begin{table}
\centering
\caption{Reduction of the total number of \ac{CG} iterations. The \ac{CG} tolerance is relaxed according to  \Cref{eq:cst,eq:z,eq:zOptim,eq:adasquare,eq:hybrid}. The parameter in \OptimalCase{} is $0.05$.  }
%\resizebox{\textwidth}{!}{%
\begin{tabular}{cccccc}
\ac{CG} tolerance  & \ConstantCase & \AdaptiveCase  & \PredictedCase  & \HybridCase & \OptimalCase \\
\ac{CG} iterations & 6963     & 5115          & 4897  & 4873  & 4399 \\
Savings \%    & -        & 26.54        & 29.67 & 30.02  & 36.82
\end{tabular}
%}
\label{tab:varTolZ}
\end{table}

\subsection{Increasing the savings by working on a simplified problem}
\label{subsec:simple}

Considering the observations in \Cref{subsec:theoPrecRelax} and the results plotted in \Cref{fig:originallower bound}, we can significantly reduce the accuracy of the inner solver only when the outer solver is in a superlinear convergence phase, when the $\abs{\zeta }$ sequence decreases rapidly. How much we can relax depends on the slope of the convergence curve. As such, to get the maximum reduction of the total number of iterations, the problem needs to be simplified, such that the convergence curve is as steep as possible and has no plateau. It is common to pair Krylov methods with other strategies, such as preconditioning, in order to improve their convergence behavior. The literature on these kinds of approaches is rich \cite{loghin2003schur,loghin2004analysis,bgl_2005,olshanskii2010acquired}. The following tests quantify how beneficial is the interaction between our proposed relaxation scheme and these other strategies.  

It has been shown by Arioli and Orban that the \ac{GKB} applied to the saddle-point system is equivalent to the \ac{CG} algorithm applied to the Schur complement equation \cite[Chapter~5]{orban2017iterative}. As such, the first step towards accelerating \ac{GKB} is to consider the Schur complement, defined as $\Sm := \Am ^T \Mm _{-1} \Am$, especially its spectrum. Ideally, a spectrum with tightly clustered values and no outliers leads to rapid \ac{GKB} convergence \cite{KrDaTaArRu2020}. To get as close as possible to this clustering we use the following two methods to induce positive changes in the spectrum: preconditioning with the \ac{LSC}  \cite{elman2006block}  and eigenvalue deflation. Each of them operates differently and leads to convergence curves with different traits. 

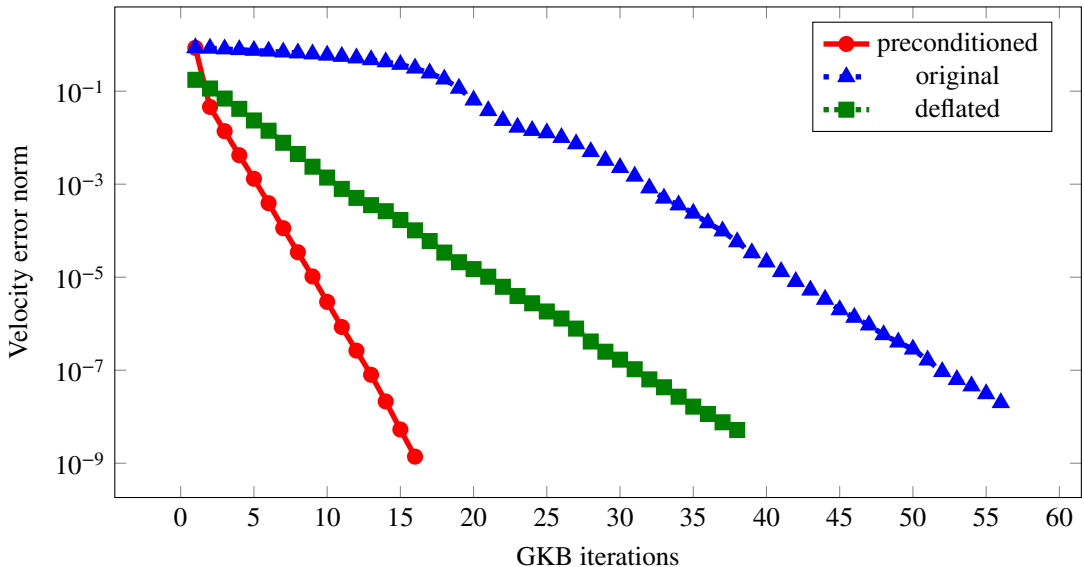
\begin{figure}
    \centering
    \begin{tikzpicture}
    \begin{axis}[
    ymode=log,
    legend pos= north east,
	table/col sep=tab,
    width=\FigWid \textwidth,
	height=\FigHei \textwidth,
    xlabel={\ac{GKB} iterations},
    ylabel={Velocity error norm},
    % xmax=70
    % extra x ticks={70},
    % extra x tick style={draw}
    ]
    
% \addplot+[] table [x=augmented LagrangianX, y=augmented LagrangianY]{images/csv/easyIFISS20GKconv.csv}; 
% \label{fig:item:easyIFISS20GKconv_augmented Lagrangian} 
% \addlegendentry{augm. Lag.} 

\addplot+[] table [x=LSC preconditionerX, y=LSC preconditionerY]{images/csv/easyIFISS20GKconv.csv}; 
\label{fig:item:easyIFISS20GKconv_LSC preconditioner} 
\addlegendentry{preconditioned} 

\addplot+[] table [x=originalX, y=originalY]{images/csv/easyIFISS20GKconv.csv}; 
\label{fig:item:easyIFISS20GKconv_original} 
\addlegendentry{original} 

\addplot+[] table [x=deflationX, y=deflationY]{images/csv/easyIFISS20GKconv.csv}; 
\label{fig:item:easyIFISS20GKconv_deflation} 
\addlegendentry{deflated} 

\end{axis} 
\end{tikzpicture} 
\caption{\ac{GKB} convergence curves for the IFISS channel test case before and after spectral clustering. Y-axis: Energy norm of the relative error for the velocity. Target \ac{GKB} tolerance \num{1e-7}. Using the Least Squares Commutator preconditioner or deflation of the smallest five spectral outliers.} 
% Augmented Lagrangian with $\eta=1000$,
\label{fig:easyIFISS20GKconv} 
\end{figure}

In \Cref{fig:easyIFISS20GKconv}, we plot the \ac{GKB} convergence curve for each of these, using a direct inner solver. The \ac{LSC} aligns the small values in the spectrum with the main cluster and brings everything closer together. The corresponding \ac{GKB} convergence curve has no plateau and is much steeper than the curve for the unpreconditioned case. Using deflation, we remove the five smallest values from the spectrum, which constitute outliers with the respect to the main cluster. The other values remain unchanged. As such, its convergence curve no longer has the initial plateau, but is otherwise the same as in the original problem. 

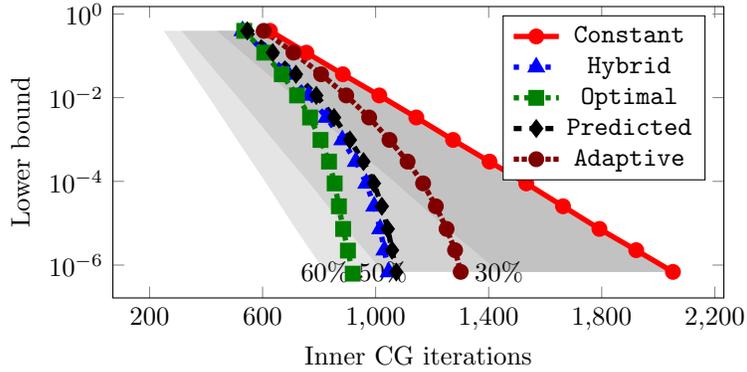
\begin{figure}
	\centering
	\begin{tikzpicture}
		\begin{axis}[
			ymode=log,
			legend pos= north east,
			table/col sep=tab,
			width=\FigWid \textwidth,
			height=\FigHei \textwidth,
			xlabel={Inner \ac{CG} iterations},
			xtick = {200,600,...,2200},
			ylabel={Lower bound}
			]
			
			\addplot+[name path=A] table [x=cumul iter, y=lower bound]{images/csv/IFISS20/LSCprec/const.csv}; 
			\addlegendentry{ \ConstantCase} 
			
			\addplot+[] table [x expr=\thisrow{cumul iter}, y=lower bound]{images/csv/IFISS20/LSCprec/hybrid.csv}; 
			\addlegendentry{\HybridCase }

			\addplot+[] table [x expr=\thisrow{cumul iter}, y=lower bound]{images/csv/IFISS20/LSCprec/zOptimal.csv}; 
			\addlegendentry{ \OptimalCase} 
			
			\addplot+[] table [x expr=\thisrow{cumul iter}, y=lower bound]{images/csv/IFISS20/LSCprec/zAdaSq.csv}; 
			\addlegendentry{ \PredictedCase} 
			
			\addplot+[] table [x expr=\thisrow{cumul iter}, y=lower bound]{images/csv/IFISS20/LSCprec/z.csv}; 
			\addlegendentry{ \AdaptiveCase} 
			
			\addplot[name path=B, draw=none] table [x expr=0.7*\thisrow{cumul iter}, y=lower bound]{images/csv/IFISS20/LSCprec/const.csv} node[pos=1]{30\%};
			
			\addplot[black!50, opacity=0.5] fill between[of=A and B];
			
			\addplot[name path=C, draw=none] table [x expr=0.5*\thisrow{cumul iter}, y=lower bound]{images/csv/IFISS20/LSCprec/const.csv} node[pos=1]{50\%};
			\addplot[black!30, opacity=0.5] fill between[of=A and C];

			\addplot[name path=C, draw=none] table [x expr=0.4*\thisrow{cumul iter}, y=lower bound]{images/csv/IFISS20/LSCprec/const.csv} node[pos=1]{60\%};
			\addplot[black!20, opacity=0.5] fill between[of=A and C];
			
		\end{axis} 
	\end{tikzpicture} 
	\caption{Lower bound (\Cref{eq:lowBnd}) for the error norm associated with the \ac{GKB} iterates versus the cumulative number of inner \ac{CG} iterations when solving the problem from \Cref{sec:pbDesc}. The problem includes preconditioning with the \ac{LSC}. The parameter used in \OptimalCase{} is $0.007$. See \Cref{eq:cst,eq:z,eq:zOptim,eq:adasquare,eq:hybrid} for the strategies denoted by the labels.}
	\label{fig:LSCpreclower bound}
\end{figure}

For both of these cases we apply the same strategies of relaxing the inner tolerance, to see how many total \ac{CG} iterations we can save. The rest of the set-up is identical to that described for \Cref{tab:varTolZ}. We tabulate the results in \Cref{tab:varTolLSC,tab:varTolDefl} and plot them in \Cref{fig:LSCpreclower bound,fig:deflatedlower bound}. They highlight that the best parameter-free results are obtained when using \HybridCase{}, which leads to savings of about 50\%, depending on the specific case. When comparing this parameter-free approach to \OptimalCase{}, which includes an experimental constant, we find that the hybrid approach can still be improved. Nonetheless, the difference in \ac{CG} iterations savings is not very high (up to 6\%), which supports the idea that our proposed strategy is efficient in a general-use setting. An additional observation pertaining to the plots is that even if convergence is relatively fast (\Cref{fig:LSCpreclower bound}) or slow (\Cref{fig:deflatedlower bound}), the final savings are still around 50\%, as long as there is no plateau.

\begin{table}
\centering
\caption{Reduction of the total number of \ac{CG} iterations after using the \ac{LSC} preconditioner. The \ac{CG} tolerance is relaxed according to  \Cref{eq:cst,eq:z,eq:zOptim,eq:adasquare,eq:hybrid}. The parameter used in \OptimalCase{} is $0.007$.  }

\begin{tabular}{cccccc}
\ac{CG} tolerance  & \ConstantCase & \AdaptiveCase  & \PredictedCase  & \HybridCase & \OptimalCase \\
\ac{CG} iterations & 2052     & 1301          & 1073  & 1046  & 919 \\
Savings \%    & -        & 36.60        & 47.71 & 49.03  & 55.21
\end{tabular}
\label{tab:varTolLSC}
\end{table}

\begin{figure}
	\centering
	\begin{tikzpicture}
		\begin{axis}[
			ymode=log,
			legend pos= north east,
			table/col sep=tab,
			width=\FigWid \textwidth,
			height=\FigHei \textwidth,
			xlabel={Inner \ac{CG} iterations},
			ylabel={Lower bound}
			]
			
			\addplot+[name path=A] table [x=cumul iter, y=lower bound]{images/csv/IFISS20/deflated/const.csv}; 
			\addlegendentry{ \ConstantCase} 
			
			\addplot+[] table [x expr=\thisrow{cumul iter}, y=lower bound]{images/csv/IFISS20/deflated/hybrid.csv}; 
			\addlegendentry{\HybridCase }

			\addplot+[] table [x expr=\thisrow{cumul iter}, y=lower bound]{images/csv/IFISS20/deflated/zOptimal.csv}; 
			\addlegendentry{ \OptimalCase} 
			
			\addplot+[] table [x expr=\thisrow{cumul iter}, y=lower bound]{images/csv/IFISS20/deflated/zAdaSq.csv}; 
			\addlegendentry{ \PredictedCase} 
			
			\addplot+[] table [x expr=\thisrow{cumul iter}, y=lower bound]{images/csv/IFISS20/deflated/z.csv}; 
			\addlegendentry{ \AdaptiveCase} 
			
			\addplot[name path=B, draw=none] table [x expr=0.6*\thisrow{cumul iter}, y=lower bound]{images/csv/IFISS20/deflated/const.csv} node[pos=1]{40\%};
			
			\addplot[black!50, opacity=0.5] fill between[of=A and B];
			
			\addplot[name path=C, draw=none] table [x expr=0.5*\thisrow{cumul iter}, y=lower bound]{images/csv/IFISS20/deflated/const.csv} node[pos=1]{50\%};
			\addplot[black!30, opacity=0.5] fill between[of=A and C];

			\addplot[name path=C, draw=none] table [x expr=0.4*\thisrow{cumul iter}, y=lower bound]{images/csv/IFISS20/deflated/const.csv} node[pos=1]{60\%};
			\addplot[black!20, opacity=0.5] fill between[of=A and C];
			
		\end{axis} 
	\end{tikzpicture} 
	\caption{Lower bound (\Cref{eq:lowBnd}) for the error norm associated with the \ac{GKB} iterates versus the cumulative number of inner \ac{CG} iterations when solving the problem from \Cref{sec:pbDesc}. The problem includes deflation of the five smallest spectral outliers. The parameter used in \OptimalCase{} is $0.09$. See \Cref{eq:cst,eq:z,eq:zOptim,eq:adasquare,eq:hybrid} for the strategies denoted by the labels.}
	\label{fig:deflatedlower bound}
\end{figure}

\begin{table}
\centering
\caption{Reduction of the total number of \ac{CG} iterations after using deflation. The \ac{CG} tolerance is relaxed according to  \Cref{eq:cst,eq:z,eq:zOptim,eq:adasquare,eq:hybrid}. 
The parameter used in \OptimalCase{} is $0.09$.  }

%\resizebox{\textwidth}{!}{%
\begin{tabular}{cccccc}
\ac{CG} tolerance  & \ConstantCase & \AdaptiveCase  & \PredictedCase  & \HybridCase & \OptimalCase \\
\ac{CG} iterations & 4830     & 2625          & 2416  & 2411  & 2110 \\
Savings \%    & -        & 45.65        & 49.98 & 50.08  & 56.31
\end{tabular}
%}
\label{tab:varTolDefl}
\end{table}

\section{\ac{GKB} with the augmented Lagrangian approach}
\label{sec:AL}
The method of the \ac{AL} has been used successfully to speed up the convergence of the \ac{GKB} algorithm \cite{KrDaTaArRu2020}, with this effect being theoretically explained by Arioli et al. \cite{KrDaTaArRu2020}.
Maybe most striking is the potential to reach mesh-independent convergence, provided that the augmentation parameter is large enough. Another use of the \ac{AL} method is to transform the (1,1)-block of a saddle-point system, say $\Wm$, from a positive semi-definite matrix to a positive definite one. However, this can happen only if the off-diagonal block $\Am$ is full rank or, more generally, if $\mbox{ker}(\Wm)\cap \mbox{ker}(\Am^T)=\{ \mZ \}$.

Let $\Nm \in \mathbb{R}^{n\times n}$ be a symmetric, positive definite matrix. For a given symmetric, positive semi-definite matrix  $\Wm \in  \mathbb{R}^{m\times m}$, we can transform it into a positive-definite one by
\begin{align}
    \Mm := \Wm +  \Am \Nm^{-1} \Am^T.
\end{align}
The upper right-hand side term $\gvv$ then becomes
\begin{equation}
\gvv := \gvv +   \Am \Nm^{-1}\rv.
\end{equation}
With these changes in place, we can proceed to using the \ac{GKB} algorithm, as described in \Cref{sec:GKBtheory}.

Note that if the matrix $\Wm$ is already symmetric positive-definite, the transformation of the (1,1)-block is not necessary for using the \ac{GKB} method. However, the application of the \ac{AL} approach does lead to a better conditioning of the Schur complement, which significantly improves convergence speed \cite{KrDaTaArRu2020}. As in \Cref{sec:GKBtheory}, we choose $\Nm=\frac{1}{\eta} \Id$. There is as usual no free lunch: depending on the conditioning of the matrix $\Am$ and the magnitude of $\eta$, the \ac{AL} can also degrade the conditioning of the $\Mm$ matrix as a side-effect. 

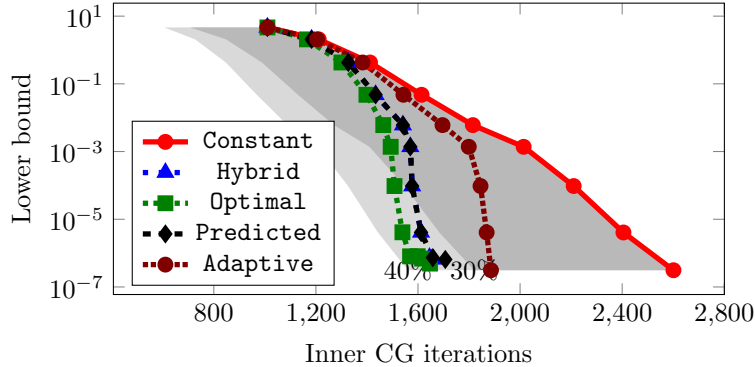
\begin{figure}
	\centering
	\begin{tikzpicture}
		\begin{axis}[
			ymode=log,
			legend pos= south west,
			table/col sep=tab,
			width=\FigWid \textwidth,
			height=\FigHei \textwidth,
			xlabel={Inner \ac{CG} iterations},
			xtick={800,1200,...,2800},
			ylabel={Lower bound}
			]
			
			\addplot+[name path=A] table [x=cumul iter, y=lower bound]{images/csv/IFISS20/augLag/const.csv}; 
			\addlegendentry{ \ConstantCase} 
			
			\addplot+[] table [x expr=\thisrow{cumul iter}, y=lower bound]{images/csv/IFISS20/augLag/hybrid.csv}; 
			\addlegendentry{\HybridCase }

			\addplot+[] table [x expr=\thisrow{cumul iter}, y=lower bound]{images/csv/IFISS20/augLag/zOptimal.csv}; 
			\addlegendentry{ \OptimalCase} 
			
			\addplot+[] table [x expr=\thisrow{cumul iter}, y=lower bound]{images/csv/IFISS20/augLag/zAdaSq.csv}; 
			\addlegendentry{ \PredictedCase} 
			
			\addplot+[] table [x expr=\thisrow{cumul iter}, y=lower bound]{images/csv/IFISS20/augLag/z.csv}; 
			\addlegendentry{ \AdaptiveCase} 
			
			\addplot[name path=B, draw=none] table [x expr=0.7*\thisrow{cumul iter}, y=lower bound]{images/csv/IFISS20/augLag/const.csv} node[pos=1]{30\%};
			
			\addplot[black!50, opacity=0.5] fill between[of=A and B];
			
			\addplot[name path=C, draw=none] table [x expr=0.6*\thisrow{cumul iter}, y=lower bound]{images/csv/IFISS20/augLag/const.csv} node[pos=1]{40\%};
			\addplot[black!30, opacity=0.5] fill between[of=A and C];

% 			\addplot[name path=C, draw=none] table [x expr=0.35*\thisrow{cumul iter}, y=lower bound]{images/csv/IFISS20/augLag/const.csv} node[pos=1]{65\%};
% 			\addplot[black!20, opacity=0.5] fill between[of=A and C];
			
		\end{axis} 
	\end{tikzpicture} 
	\caption{Lower bound (\Cref{eq:lowBnd}) for the error norm associated with the \ac{GKB} iterates versus the cumulative number of inner \ac{CG} iterations when solving the problem from \Cref{sec:pbDesc}. The problem includes the \ac{AL} ($\eta =1000$). The parameter used in \OptimalCase{} is $0.005$. See \Cref{eq:cst,eq:z,eq:zOptim,eq:adasquare,eq:hybrid} for the strategies denoted by the labels.}
	\label{fig:augLaglower bound}
\end{figure}

We test whether the augmentation interacts with the strategies we propose in \Cref{sec:relaxChoices}, namely if we can still achieve about 50\% savings in the total number of inner iterations. The strategies are applied when solving the problem described in \Cref{sec:pbDesc} after an augmentation with a parameter $\eta=1000$, with the results being given in \Cref{tab:varTolAL} and plotted in \Cref{fig:augLaglower bound}. Comparing the percentage of iterations saved in this case to those obtained in \Cref{sec:relaxChoices}, it is clear that, when combined with the \ac{AL} method, the strategy of variable inner tolerance does help reducing the total number of inner iterations, but by a lower percentage. 

\begin{table}
\centering
\caption{Reduction of the total number of \ac{CG} iterations after using the \ac{AL} ($\eta =1000$). The \ac{CG} tolerance is relaxed according to  \Cref{eq:cst,eq:z,eq:zOptim,eq:adasquare,eq:hybrid}. 
The parameter used in \OptimalCase{} is $0.005$.  }
%\resizebox{\textwidth}{!}{%
\begin{tabular}{cccccc}
\ac{CG} tolerance  & \ConstantCase & \AdaptiveCase  & \PredictedCase  & \HybridCase & \OptimalCase \\
\ac{CG} iterations & 2601     & 1886          & 1707  & 1661  & 1647 \\
Savings \%    & -        & 27.49        & 34.37 & 36.14 & 36.68
\end{tabular}
%}
\label{tab:varTolAL}
\end{table}

Since the \ac{AL} method modifies the (1,1)-block of the saddle-point system, it changes the difficulty of the inner problem and how many iterations the inner solver needs to perform. As such, a global comparison in terms of number of inner iterations, among all the scenarios we studied (original, preconditioned, deflated, including the \ac{AL}) is not fair unless the inner problem has the same degree of difficulty for all the cases.

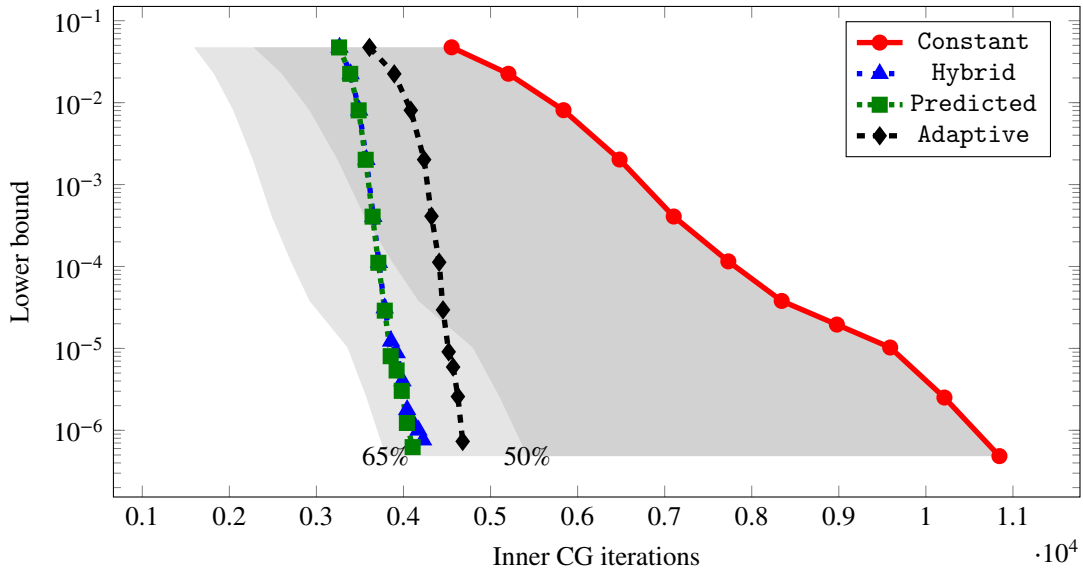
\begin{figure}
	\centering
	\begin{tikzpicture}
		\begin{axis}[
			ymode=log,
			legend pos= north east,
			%table/col sep=tab,
			width=\FigWid \textwidth,
			height=\FigHei \textwidth,
			xlabel={Inner CG iterations},
			ylabel={Lower bound}
			]
			
			\addplot+[name path=A] table [x=InnerKSPCumul, y=lowbnd]{images/petscResults/MixedPoisson/monitorCGJacobiCST500.txt}; 
% 			\addlegendentry{Constant accuracy} 
			\addlegendentry{ \ConstantCase}

			\addplot+[] table [x expr=\thisrow{InnerKSPCumul}, y=lowbnd]{images/petscResults/MixedPoisson/monitorCGJacobiHybrid500.txt}; 
% 			\addlegendentry{hybrid accuracy} 
			\addlegendentry{\HybridCase }

% 			\addplot+[] table [x expr=\thisrow{InnerKSPCumul}, y=lowbnd]{images/petscResults/MixedPoisson/monitorCGJacobiAda500.txt}; 
% 			\addlegendentry{Ada accuracy} 
			
			\addplot+[] table [x expr=\thisrow{InnerKSPCumul}, y=lowbnd]{images/petscResults/MixedPoisson/monitorCGJacobiAdaSquare500.txt}; 
			\addlegendentry{ \PredictedCase} 
% 			\addlegendentry{Ada Squared accuracy} 
			
			\addplot+[] table [x expr=\thisrow{InnerKSPCumul}, y=lowbnd]{images/petscResults/MixedPoisson/monitorCGJacobiZ500.txt}; 
% 			\addlegendentry{z accuracy} 
			\addlegendentry{ \AdaptiveCase}

% 			\addplot[name path=B, draw=none] table [x expr=0.7*\thisrow{InnerKSPCumul}, y=lowbnd]{images/petscResults/MixedPoisson/monitorCGJacobiCST500.txt} node[pos=1]{30\%};
% 			\addplot[black!50, opacity=0.5] fill between[of=A and B];
			
			\addplot[name path=C, draw=none] table [x expr=0.5*\thisrow{InnerKSPCumul}, y=lowbnd]{images/petscResults/MixedPoisson/monitorCGJacobiCST500.txt} node[pos=1]{50\%};
			\addplot[black!30, opacity=0.5] fill between[of=A and C];

			\addplot[name path=C, draw=none] table [x expr=0.35*\thisrow{InnerKSPCumul}, y=lowbnd]{images/petscResults/MixedPoisson/monitorCGJacobiCST500.txt} node[pos=1]{65\%};
			\addplot[black!20, opacity=0.5] fill between[of=A and C];
			
		\end{axis} 
	\end{tikzpicture} 
	
	\caption{Lower bound (\Cref{eq:lowBnd}) for the error norm associated with the \ac{GKB} iterates versus the cumulative number of inner \ac{CG} iterations when solving the Mixed Poisson problem. We also use the \ac{AL} ($\eta =500$). See \Cref{eq:cst,eq:z,eq:adasquare,eq:hybrid} for the strategies denoted by the labels.}
	
	\label{fig:mixedPoissonLowBnd}
\end{figure}

To verify the generality of our method, we also apply it in a different context than that described in \Cref{sec:pbDesc}. Let us consider a Mixed Poisson problem. We solve the Poisson equation \(-\Delta u=f\) on the unit square $(0,1)^2$ using a mixed formulation. We introduce the vector variable $\vec{\sigma}=\nabla u$. Find $(\vec{\sigma}, u)\in \Sigma \times W $ such that 
\begin{align}
    \vec{\sigma}-\grad u&=0\\
    -\mathrm{div} (\vec{\sigma}) &= f.
\end{align}
where homogeneous Dirichlet boundary conditions are imposed for $u$ at all walls. The forcing term $f$ is random and uniformly drawn in $(0,1)$. The discretization is done with a lowest order Raviart-Thomas space $\Sigma^h \subset \Sigma$, and a space $W^h \subset W$ containing piece-wise constant basis functions. We used the finite element package Firedrake\footnote{\url{www.firedrakeproject.org}} coupled with a PETSc~\cite{petsc-web-page,petsc-user-ref,petsc-efficient} implementation of \ac{GKB} \footnote{\url{https://petsc.org/release/docs/manualpages/PC/PCFIELDSPLIT.html\#PCFIELDSPLIT}}, adapted to include dynamical relaxation, to produce the following numerical results. We used the implementation provided by Firedrake\footnote{\url{https://www.firedrakeproject.org/demos/saddle_point_systems.py.html}}.
The test case has \num{328192} degrees of freedom, of which \num{197120} are associated with the (1,1)-block. The \ac{GKB} delay parameter is set to 3. The augmentation parameter $\eta$ is set to 500 and the tolerance for the \ac{GKB} set to \num{1e-5}. The results are presented in \Cref{fig:mixedPoissonLowBnd}. We confirm the results presented above with a reduction of over 60\% in the total number of inner \ac{CG} iterations with respect to the constant accuracy set up.

\begin{table}
\centering
\caption{Reduction of the total number of \ac{CG} iterations after using the \ac{AL} ($\eta =500$) on the Mixed Poisson problem. The \ac{CG} tolerance is relaxed according to  \Cref{eq:cst,eq:z,eq:adasquare,eq:hybrid}. }
%\resizebox{\textwidth}{!}{%
\begin{tabular}{ccccc}
\ac{CG} tolerance  & \ConstantCase & \AdaptiveCase  & \PredictedCase  & \HybridCase  \\
\ac{CG} iterations & 10845     & 4680          & 4105  & 4225    \\
Savings \%    & -        & 56.84        & 62.15 & 61.04  
\end{tabular}
%}
\label{tab:DarcyvarTolAL}
\end{table}

\section{Conclusions}

We have studied the behavior of the \ac{GKB} algorithm in the case where the inner problem, i.e. the solution of a linear system, is performed iteratively. 
We have found that the inner solver does not need to be as precise as a direct one in order to achieve a \ac{GKB} solution of a predefined accuracy. 

Furthermore, we have proposed algorithmic strategies that reduce the cost of the inner solver, quantified as the cumulative number of inner iterations.
This is possible by selecting criteria to change the stopping tolerance.
To motivate these choices, we have studied the perturbation generated by the inexact inner solver. The findings show that the perturbation introduced in early iterations has a higher impact on the accuracy of the solution compared to later ones.

We devised a dynamic way of adapting the accuracy of the inner solver at each call to minimize its cost. The initial, high accuracy is gradually reduced, maintaining the resulting perturbation under control. 

Our relaxation strategy is inexpensive, easy to implement, and has reduced the total number of inner iterations by 33-63\% in our tests. The experiments also show that including methods such as deflation, preconditioning and the augmented Lagrangian has no negative impact and can lead to a higher percentage of savings. Another advantage is that our method does not rely on additional parameters and is thus usable in a black-box fashion.

\paragraph{Acknowledgments}
The authors thank Mario Arioli for many inspiring discussions and advice.

\bibliography{literature}

\begin{thebibliography}{10}

\bibitem{Ar2013}
M.~Arioli.
\newblock Generalized {G}olub--{K}ahan bidiagonalization and stopping criteria.
\newblock {\em SIAM Journal on Matrix Analysis and Applications},
  34(2):571--592, 2013.

\bibitem{petsc-user-ref}
Satish Balay, Shrirang Abhyankar, Mark~F. Adams, Steven Benson, Jed Brown,
  Peter Brune, Kris Buschelman, Emil Constantinescu, Lisandro Dalcin, Alp
  Dener, Victor Eijkhout, William~D. Gropp, V\'{a}clav Hapla, Tobin Isaac,
  Pierre Jolivet, Dmitry Karpeev, Dinesh Kaushik, Matthew~G. Knepley, Fande
  Kong, Scott Kruger, Dave~A. May, Lois~Curfman McInnes, Richard~Tran Mills,
  Lawrence Mitchell, Todd Munson, Jose~E. Roman, Karl Rupp, Patrick Sanan,
  Jason Sarich, Barry~F. Smith, Stefano Zampini, Hong Zhang, Hong Zhang, and
  Junchao Zhang.
\newblock {PETSc/TAO} users manual.
\newblock Technical Report ANL-21/39 - Revision 3.16, Argonne National
  Laboratory, 2021.

\bibitem{petsc-web-page}
Satish Balay, Shrirang Abhyankar, Mark~F. Adams, Steven Benson, Jed Brown,
  Peter Brune, Kris Buschelman, Emil~M. Constantinescu, Lisandro Dalcin, Alp
  Dener, Victor Eijkhout, William~D. Gropp, V\'{a}clav Hapla, Tobin Isaac,
  Pierre Jolivet, Dmitry Karpeev, Dinesh Kaushik, Matthew~G. Knepley, Fande
  Kong, Scott Kruger, Dave~A. May, Lois~Curfman McInnes, Richard~Tran Mills,
  Lawrence Mitchell, Todd Munson, Jose~E. Roman, Karl Rupp, Patrick Sanan,
  Jason Sarich, Barry~F. Smith, Stefano Zampini, Hong Zhang, Hong Zhang, and
  Junchao Zhang.
\newblock {PETS}c {W}eb page.
\newblock \url{https://petsc.org/}, 2021.

\bibitem{petsc-efficient}
Satish Balay, William~D. Gropp, Lois~Curfman McInnes, and Barry~F. Smith.
\newblock Efficient management of parallelism in object oriented numerical
  software libraries.
\newblock In E.~Arge, A.~M. Bruaset, and H.~P. Langtangen, editors, {\em Modern
  Software Tools in Scientific Computing}, pages 163--202. Birkh{\"{a}}user
  Press, 1997.

\bibitem{baumann2015nested}
Manuel Baumann and Martin~B Van~Gijzen.
\newblock Nested krylov methods for shifted linear systems.
\newblock {\em SIAM Journal on Scientific Computing}, 37(5):S90--S112, 2015.

\bibitem{bgl_2005}
Michele Benzi, Gene~H. Golub, and J{\"o}rg Liesen.
\newblock Numerical solution of saddle point problems.
\newblock {\em Acta Numerica}, 14:1--137, 2005.

\bibitem{bouras2005inexact}
Amina Bouras and Val{\'e}rie Frayss{\'e}.
\newblock Inexact matrix-vector products in krylov methods for solving linear
  systems: a relaxation strategy.
\newblock {\em SIAM Journal on Matrix Analysis and Applications},
  26(3):660--678, 2005.

\bibitem{bouras2000relaxation}
Amina Bouras, Val{\'e}rie Frayss{\'e}, and Luc Giraud.
\newblock A relaxation strategy for inner-outer linear solvers in domain
  decomposition methods, 2000.
\newblock Technical Report 17.

\bibitem{chung2019flexible}
Julianne Chung and Silvia Gazzola.
\newblock Flexible krylov methods for $\ell_p$ regularization.
\newblock {\em SIAM Journal on Scientific Computing}, 41(5):S149--S171, 2019.

\bibitem{dax2019restarted}
Achiya Dax.
\newblock A restarted krylov method with inexact inversions.
\newblock {\em Numerical Linear Algebra with Applications}, 26(1):e2213, 2019.

\bibitem{dembo1982inexact}
Ron~S Dembo, Stanley~C Eisenstat, and Trond Steihaug.
\newblock Inexact {N}ewton methods.
\newblock {\em SIAM Journal on Numerical Analysis}, 19(2):400--408, 1982.

\bibitem{elman2006block}
Howard Elman, Victoria~E Howle, John Shadid, Robert Shuttleworth, and Ray
  Tuminaro.
\newblock Block preconditioners based on approximate commutators.
\newblock {\em SIAM Journal on Scientific Computing}, 27(5):1651--1668, 2006.

\bibitem{ers07}
Howard Elman, Alison Ramage, and David Silvester.
\newblock Algorithm {866}: {IFISS}, a {M}atlab toolbox for modelling
  incompressible flow.
\newblock {\em ACM Trans. Math. Softw.}, 33:2--14, 2007.

\bibitem{elman2014finite}
Howard~C Elman, David~J Silvester, and Andrew~J Wathen.
\newblock {\em {Finite elements and fast iterative solvers: with applications
  in incompressible fluid dynamics}}.
\newblock Numerical Mathematics and Scie, 2014.

\bibitem{erlangga2008multilevel}
Yogi~A Erlangga and Reinhard Nabben.
\newblock Multilevel projection-based nested krylov iteration for boundary
  value problems.
\newblock {\em SIAM Journal on Scientific Computing}, 30(3):1572--1595, 2008.

\bibitem{gazzola2021regularization}
Silvia Gazzola and Malena~Sabate Landman.
\newblock Regularization by inexact krylov methods with applications to blind
  deblurring.
\newblock {\em SIAM Journal on Matrix Analysis and Applications},
  42(4):1528--1552, 2021.

\bibitem{GolZhaZha2000}
Gene~H. Golub, Zhenyue Zhang, and Hongyuan Zha.
\newblock Large sparse symmetric eigenvalue problems with homogeneous linear
  constraints: the lanczos process with inner–outer iterations.
\newblock {\em Linear Algebra and its Applications}, 309(1):289 -- 306, 2000.

\bibitem{jiranek2008maximum}
Pavel Jir{\'a}nek and Miroslav Rozlo{\v{z}}n{\'\i}k.
\newblock Maximum attainable accuracy of inexact saddle point solvers.
\newblock {\em SIAM journal on matrix analysis and applications},
  29(4):1297--1321, 2008.

\bibitem{kehl2019adaptive}
Ren{\'e} Kehl, Reinhard Nabben, and Daniel~B Szyld.
\newblock Adaptive multilevel krylov methods.
\newblock {\em Electronic Transactions on Numerical Analysis}, 51, 2019.

\bibitem{KrDaTaArRu2020}
C.~Kruse, V.~Darrigrand, N.~Tardieu, M.~Arioli, and U.~Rüde.
\newblock Application of an iterative golub-kahan algorithm to structural
  mechanics problems with multi-point constraints.
\newblock {\em Adv. Model. and Simul. in Eng. Sci}, 7, 2020.

\bibitem{loghin2003schur}
Daniel Loghin and Andrew~J Wathen.
\newblock Schur complement preconditioning for elliptic systems of partial
  differential equations.
\newblock {\em Numerical linear algebra with applications}, 10(5-6):423--443,
  2003.

\bibitem{loghin2004analysis}
Daniel Loghin and Andrew~J Wathen.
\newblock Analysis of preconditioners for saddle-point problems.
\newblock {\em SIAM Journal on Scientific Computing}, 25(6):2029--2049, 2004.

\bibitem{mcinnes2014hierarchical}
Lois~Curfman McInnes, Barry Smith, Hong Zhang, and Richard~Tran Mills.
\newblock Hierarchical krylov and nested krylov methods for extreme-scale
  computing.
\newblock {\em Parallel Computing}, 40(1):17--31, 2014.

\bibitem{olshanskii2010acquired}
Maxim~A Olshanskii and Valeria Simoncini.
\newblock Acquired clustering properties and solution of certain saddle point
  systems.
\newblock {\em SIAM journal on matrix analysis and applications},
  31(5):2754--2768, 2010.

\bibitem{orban2017iterative}
Dominique Orban and Mario Arioli.
\newblock {\em {Iterative solution of symmetric quasi-definite linear
  systems}}.
\newblock SIAM, 2017.

\bibitem{saad1993flexible}
Youcef Saad.
\newblock A flexible inner-outer preconditioned gmres algorithm.
\newblock {\em SIAM Journal on Scientific Computing}, 14(2):461--469, 1993.

\bibitem{simoncini2003theory}
Valeria Simoncini and Daniel~B Szyld.
\newblock Theory of inexact krylov subspace methods and applications to
  scientific computing.
\newblock {\em SIAM Journal on Scientific Computing}, 25(2):454--477, 2003.

\bibitem{simoncini2005relaxed}
Valeria Simoncini and Daniel~B Szyld.
\newblock Relaxed krylov subspace approximation.
\newblock In {\em PAMM: Proceedings in Applied Mathematics and Mechanics},
  volume~5, pages 797--800. Wiley Online Library, 2005.

\bibitem{simoncini2007recent}
Valeria Simoncini and Daniel~B Szyld.
\newblock Recent computational developments in krylov subspace methods for
  linear systems.
\newblock {\em Numerical Linear Algebra with Applications}, 14(1):1--59, 2007.

\bibitem{van2004inexact}
Jasper Van Den~Eshof and Gerard~LG Sleijpen.
\newblock Inexact krylov subspace methods for linear systems.
\newblock {\em SIAM Journal on Matrix Analysis and Applications},
  26(1):125--153, 2004.

\bibitem{van2005relaxation}
Jasper Van Den~Eshof, Gerard~LG Sleijpen, and Martin~B van Gijzen.
\newblock Relaxation strategies for nested krylov methods.
\newblock {\em Journal of Computational and Applied Mathematics},
  177(2):347--365, 2005.

\bibitem{van1994gmresr}
Henk~A Van~der Vorst and Cornelis Vuik.
\newblock Gmresr: a family of nested gmres methods.
\newblock {\em Numerical Linear Algebra with Applications}, 1(4):369--386,
  1994.

\bibitem{xu2022inexact}
Shengjie Xu and Fei Xue.
\newblock Inexact rational krylov subspace method for eigenvalue problems.
\newblock {\em Numerical Linear Algebra with Applications}, page e2437, 2022.

\end{thebibliography}
\bibliographystyle{plain}

\end{document}